\documentclass[12pt]{article}

\usepackage[affil-it]{authblk}
      \usepackage{latexsym}
         \usepackage[reqno, namelimits, sumlimits]{amsmath}
         \usepackage{amssymb, amsfonts}
         \usepackage{amsthm}

 \newtheorem{theorem}{Theorem}[section]

\title{Liouville Type Theorem for Stationary Navier-Stokes Equations}
\author{ G Seregin
  }
\affil{OxPDE, Mathematical Institute, University of Oxford, Oxford,UK}
\date{ \today}

\begin{document}
\maketitle
\begin{abstract}
It is shown that any smooth solution to the stationary Navier-Stokes system in $R^3$ with the velocity field, belonging globally to $L_6$ and $BM0^{-1}$, must be zero.\end{abstract}\


\setcounter{equation}{0}
\section{Introduction}

Let us consider smooth solutions to the  stationary Navier-Stokes system
\begin{equation}\label{nse}
u\cdot \nabla u-\Delta u=\nabla p, \qquad {\rm div }u=0
\end{equation}
in $\mathbb R^3$ with the additional condition at infinity: $u(x)\to 0$ as $|x|\to \infty$. Then the question is weather or not $u$ must be identically zero.
The point addressed in this short note is under which additional assumptions the answer to the above question is positive. 

In the monograph \cite{Galdi-book}, it has been shown that the condition 
\begin{equation}\label{9/2}
u\in L_\frac 92(\mathbb R^3)
\end{equation}
implies $u=0$.

A plausible conjecture is that a sufficient condition for the positive answer could be as follows
\begin{equation}\label{finite dissipation}
\int\limits_{\mathbb R^3}|\nabla u|^2dx<\infty.
\end{equation}

At the moment of writing this note, a proof of whether or not   (\ref{finite dissipation})
is a sufficient condition for $u$  to be identically zero is not known yet. We are going to show that however under an additional condition it is true. To describe that condition, we need the following definition.
We say that a divergence free vector valued field $u$ belongs to the space $BMO^{-1}(\mathbb R^3)$
 if there exists skew symmetric tensor $d$ from $BMO(\mathbb R^3)$ such that $$
 u={\rm div}\,d=(d_{ij,j}).$$
 
  It is known, see for instance \cite{Stein1970}, that if $d\in BMO(\mathbb R^3)$, then 
  $$  \Gamma(s):=\sup\limits_{x_0\in\mathbb R^3,0<r}\Big(\frac 1{|B(r)|}\int\limits_{B(x_0,r)}|d-[d]_{x_0,r}|^sdx \Big)^\frac 1s<\infty$$
for each $1\leq s<\infty$. Here,  we denote by $B(x_0,r)$ the ball of radius $r$ centred at point $x_0$ and
$[d]_{x_0,r}$ is the mean value of $d$ over the ball $B(x_0,r)$.

Our aim is to prove the following result.
\begin{theorem}\label{main result}
The following  statements are true:

\noindent
(A) any smooth divergence free vector-valued field $u\in BMO^{-1}(\mathbb R^3)$, satisfying condition (\ref{finite dissipation}) and the system (\ref{nse}), is identically equal to zero; 


\noindent
(B) 
any smooth  divergence free vector-valued field $u\in BMO^{-1}(\mathbb R^3)$, satisfying the condition
\begin{equation}\label{6}
u\in L_6(\mathbb R^3)
\end{equation}
and the system (\ref{nse}), is identically equal to zero.
\end{theorem}

By the known inequality
$$\|v\|_{6,\Omega}\leq c\|\nabla v\|_{2,\Omega},$$
being valid for any $v\in C^\infty_0(\mathbb R^3 ) $, the statement (A) follows from  the statement (B).

Before proving Theorem \ref{main result}, let us show that the assumptions in  (B) do not follow immediately from the condition (\ref{9/2}  )  of \cite{Galdi-book}. Indeed, we can let 
$$w=\sin \Big((|x|^2+1)^{\frac 14-\varepsilon }\Big)(1,1,1)$$
and 
$$v= {\rm rot} \,w.$$
Direct calculations show that $v\in BMO^{-1}(\mathbb R^3)\cap L_6(\mathbb R^3 )   $ but
$v \notin L_\frac 92(\mathbb R^3 ) $.
 
 Finally, we would like to mention that there are many interesting papers devoted to the above or related questions, see, for example, \cite{GilWein1978}, \cite{KNSS2009}, \cite{ChaeYoneda2013}, and \cite{Chae2014}.

\setcounter{equation}{0}
\section{Proof of Main Result}
\subsection{Caccioppoli Type Inequality}
This is the main technical part of the proof. We take an arbitrary ball $B(x_0,R)\in \mathbb R^3$ and a non-negative cut off function $\varphi\in C^\infty_0(B(x_0,R))$ with the following properties  $\varphi(x)=1$ in $B(x_0,\varrho  )$, $\varphi(x)=0$ out of $B(x_0,r)$, and $|\nabla \varphi(x)|\leq c/(r-\varrho)$ for any $R/2\leq\varrho<r\leq R$. 

We let 
$$  \overline{u}=u-u_0,\qquad \overline{d}=d-[d]_{x_0,R},$$
where $u_0$  is an arbitrary constant.

We also know that, for any $2<s<\infty$, there exists a  constant $c_0=c_0(s)>0$ and a function $w\in W^1_s(B(x_0,r))$, vanishing on $\partial B(x_0,r)$,  such that ${\rm div}\,w=\nabla \varphi\cdot \overline u$ and 
\begin{equation}\label{Bogovskii}
\int\limits_{B(x_0,r)}|\nabla w|^sdx\leq c_0\int\limits_{B(x_0,r)}|\nabla \varphi\cdot \overline u|^sdx\leq \frac {c_0} {(r-\varrho)^s}\int\limits_{B(x_0,R)}| \overline u|^sdx.\end{equation}

Now, we can test the Navier-Stokes equations (\ref{nse})  with the function $\varphi \overline u-w$, integrate by parts in $B(x_0,r)$, and find the following identity
$$\int\limits_{B(x_0,r)}\varphi |\nabla u|^2dx=
-\int\limits_{B(x_0,r)}\nabla u :(\nabla \varphi\otimes\overline u) dx+\int\limits_{B(x_0,r)}\nabla w :\nabla u dx+$$
$$-\int\limits_{B(x_0,r)}(u\cdot\nabla u)\cdot\varphi\overline u dx+\int\limits_{B(x_0,r)}(u\cdot\nabla u)\cdot wdx=I_1+I_2+I_3+I_4.$$

$I_1$ and $I_2$ can be estimated easily. As a result, we find
$$|I_1|+|I_2|\leq \Big(\int\limits_{B(x_0,r)}|\nabla u|^2dx\Big)^\frac 12\frac{R^{3\frac {s-2}{2s}}}{r-\varrho}\Big(\int\limits_{B(x_0,R)}|\overline u|^sdx\Big)^\frac 1s. $$

To estimate $I_3$ and $I_4$, we are going to use the skew symmetry of the matrix $d$. We have
$$|I_3|=\Big| \int\limits_{B(x_0,r)} \overline d_{jm,m}\overline u_{i,j}\overline u_i\varphi dx\Big|=\Big| \int\limits_{B(x_0,r)} \overline d_{jm}\overline u_{i,j}\overline u_i\varphi_{,m} dx\Big|\leq $$
$$\leq \frac 1{r-\varrho}\Big(\int\limits_{B(x_0,r)}|\nabla u|^2dx\Big)^\frac 12
\Big(\int\limits_{B(x_0,r)}|\overline d|^2|\overline u|^2dx\Big)^\frac 12\leq $$
$$\leq \frac 1{r-\varrho}\Big(\int\limits_{B(x_0,r)}|\nabla u|^2dx\Big)^\frac 12
\Big(\int\limits_{B(x_0,r)}|\overline u|^sdx\Big)^\frac 1s\Big(\int\limits_{B(x_0,r)}|\overline d|^\frac {2s}{s-2}dx\Big)^\frac {s-2}{2s}\leq $$
$$\leq c\frac {R^{3\frac {s-2}{2s}}}{r-\varrho}\Big(\int\limits_{B(x_0,r)}|\nabla u|^2dx\Big)^\frac 12
\Big(\int\limits_{B(x_0,R)}|\overline u|^sdx\Big)^\frac 1s\Gamma(2s/(s-2)).$$
It remains to evaluated $I_4$:
$$|I_4|=\Big| \int\limits_{B(x_0,r)} \overline d_{jm,m}\overline u_{i,j} w_idx\Big|=\Big| \int\limits_{B(x_0,r)} \overline d_{jm}\overline u_{i,j} w_{i,m} dx\Big|\leq $$
$$\leq \frac 1{r-\varrho}\Big(\int\limits_{B(x_0,r)}|\nabla u|^2dx\Big)^\frac 12
\Big(\int\limits_{B(x_0,r)}|\overline d|^2|\nabla w|^2dx\Big)^\frac 12\leq $$
$$\leq \frac 1{r-\varrho}\Big(\int\limits_{B(x_0,r)}|\nabla u|^2dx\Big)^\frac 12
\Big(\int\limits_{B(x_0,r)}|\nabla w|^sdx\Big)^\frac 1s\Big(\int\limits_{B(x_0,r)}|\overline d|^\frac {2s}{s-2}dx\Big)^\frac {s-2}{2s}\leq $$
$$\leq c(s)\frac {R^{3\frac {s-2}{2s}}}{r-\varrho}\Big(\int\limits_{B(x_0,r)}|\nabla u|^2dx\Big)^\frac 12
\Big(\int\limits_{B(x_0,R)}|\overline u|^sdx\Big)^\frac 1s\Gamma(2s/(s-2)).$$
So, we have  the following inequality
$$\int\limits_{B(x_0,\varrho )}|\nabla u|^2dx
\leq c(s)\frac {R^{3\frac {s-2}{2s}}}{r-\varrho}\Big(\int\limits_{B(x_0,r)}|\nabla u|^2dx\Big)^\frac 12
\Big(\int\limits_{B(x_0,R)}|\overline u|^sdx\Big)^\frac 1s.$$
We can apply Young's inequality and conclude
$$\int\limits_{B(x_0,\varrho )}|\nabla u|^2dx\leq \frac 14
\int\limits_{B(x_0,r)}|\nabla u|^2dx
+c(s)\frac 1 {(r-\varrho)^2}R^{3\frac {s-2}{s}}
\Big(\int\limits_{B(x_0,R)}|\overline u|^sdx\Big)^\frac 2s.$$
The later inequality is valid for any $R/2\leq \varrho<r\leq R$.  It is known,  see for instance \cite{Giaquinta1983} (this is just a matter of  suitable iterations), that such an inequality implies the following Caccioppoli type inequality
\begin{equation}\label{Caccioppoli}
\int\limits_{B(x_0,R/2 )}|\nabla u|^2dx\leq c(s)R^{3\frac {s-2}{s}-2}
\Big(\int\limits_{B(x_0,R)}|\overline u|^sdx\Big)^\frac 2s
\end{equation}
which holds for any $B(x_0,R)$ in $\mathbb R^3$.

\subsection{(A) implies (B)}
Indeed, if we let $s=6$ and $u_0=0$, then (\ref{Caccioppoli}) takes the form
$$\int\limits_{B(x_0,R/2 )}|\nabla u|^2dx\leq c(s)\Big(\int\limits_{B(x_0,R)}| u|^6dx\Big)^\frac 13\leq c(s)\Big(\int\limits_{\mathbb R^3}| u|^6dx\Big)^\frac 13.$$
Passing $R\to \infty$, we conclude that (\ref{finite dissipation}) holds.
\subsection{Proof of (A)}
Now, we let $s=3$ and $u_0=[u]_{x_0,R}$ and use  the Gagliardo-Nirenberg type inequality
$$
\Big(\int\limits_{B(x_0,R)}|\overline u|^3dx\Big)^\frac 13\leq c\Big(\int\limits_{B(x_0,R )}|\nabla u|^\frac 32dx\Big)^\frac 23 $$ with a universal positive constant $c$. Now, (\ref{Caccioppoli}) can reduced to the following reverse H\"older inequality
$$\frac 1{|B(R/2)|}\int\limits_{B(x_0,R/2 )}|\nabla u|^2dx\leq c\Big(\frac 1{|B(R)|}\int\limits_{B(x_0,R )}|\nabla u|^\frac 32dx\Big)^\frac 43 $$
with a constant $c$ that is independent of $x_0$ and $R$. We let $h:=|\nabla u|^\frac 32\in L_\frac 43(\mathbb R^3)$ and let $M(h)$ is the maximal function of $h$, i.e.,
$$M_h(x_0)=\sup_{R>0}\int\limits_{B(x_0,R )}h(x)dx.$$
Then from the above inequality, it follows that
$$M_{h^\frac 43}(x_0)\leq cM^\frac 43_h(x_0).$$
It is known, see \cite{Stein1970}, that 
the right hand side of the latter inequality is integrable in $\mathbb R^3$ and the corresponding integral is bounded from above by the quantity
$$\int\limits_{B(x_0,R )}h^\frac 43dx=
\int\limits_{B(x_0,R )}|\nabla u|^2dx$$
times a universal constant.
So, this means that $M_{h^\frac 43}\in L_1(\mathbb R^3)$. Since $h^\frac 43\in L_1(\mathbb R^3)$, it is possible only if $h$ is identically equal to zero, see \cite{Stein1970}. So, $\nabla u$ is identically equal to zero.

\end{document}